\documentclass[12pt,oneside]{amsart}
\usepackage[dvipsnames]{xcolor}

\usepackage{comment,enumitem}

\usepackage{todonotes}

\newcommand{\philt}[1]{\todo[inline, color=blue!30]{Philippe\\ #1}}

\usepackage{mathtools}
\mathtoolsset{showonlyrefs}
\usepackage{bbm}

\numberwithin{equation}{section}

\usepackage[margin=1in]{geometry}

\usepackage{graphicx}
\usepackage{amssymb}
\usepackage{tikz}
\usepackage{amsmath,amsthm,amscd,amssymb}
\usepackage{latexsym}
\usepackage[colorlinks,citecolor=red,pagebackref,hypertexnames=false]{hyperref}
\hypersetup{%
  colorlinks = true,
  citecolor=red,
  linkcolor  = black
}

\newcommand{\bA}{\mathbb{A}}

\newcommand{\bC}{\mathbb{C}}

\newcommand{\bR}{\mathbb{R}}
\newcommand{\bS}{\mathbb{S}}
\newcommand{\bZ}{\mathbb{Z}}
\newcommand{\bT}{\mathbb{T}}

\newcommand{\cF}{\mathcal{F}}

\newcommand{\cS}{\mathcal{S}}

\newcommand{\norm}[1]{{\left\Vert#1\right\Vert}}

\newcommand{\scal}[1]{{\left\langle{#1}\right\rangle}}

\DeclareMathOperator{\supp}{supp}

\newcommand{\eps}{{\varepsilon}}

\usepackage{soul}
\numberwithin{equation}{section}
\theoremstyle{plain}
\newtheorem{theorem}{Theorem}[section]

\newtheorem{corollary}[theorem]{Corollary}
\newtheorem{proposition}[theorem]{Proposition}
\newtheorem{conjecture}[theorem]{Conjecture}

\theoremstyle{definition}
\newtheorem{definition}[theorem]{Definition}

\theoremstyle{remark}
\newtheorem{remark}[theorem]{Remark}

\newtheorem{case[theorem]}{Case}

\def\R{\mathbb R}

\date{\today}      
\author{A. Iosevich, P. Jaming and A. Mayeli} 
\address{Department of Mathematics, University of Rochester, Rochester, NY}
\email{iosevich@gmail.com}
\address{Univ. Bordeaux, CNRS, Bordeaux INP, IMB, UMR 5251, F-33400 Talence, France}
\email{philippe.jaming@math.u-bordeaux.fr}
\address{Department of Mathematics, City University of New York, NY} 
\email{amayeli@gc.cuny.edu} 

\begin{document}
\nocite{*}

\title[Annihilating pairs]{Uncertainty Principle, annihilating pairs and Fourier restriction}

\begin{abstract}
Let \(G\) be a locally compact abelian group, and let \(\widehat{G}\) denote its dual group, equipped with a Haar measure.  
A variant of the uncertainty principle states that for any \(S \subset G\) and \(\Sigma \subset \widehat{G}\), there exists a constant \(C(S, \Sigma)\) such that for any \(f \in L^2(G)\), the following inequality holds:  
\[
\|f\|_{L^2(G)} \leq C(S, \Sigma) \bigl( \|f\|_{L^2(G \setminus S)} + \|\widehat{f}\|_{L^2(\widehat{G} \setminus \Sigma)} \bigr),
\]  
where \(\widehat{f}\) denotes the Fourier transform of \(f\).  
This variant of the uncertainty principle is particularly useful in applications such as signal processing and control theory.

The purpose of this paper is to show that such estimates can be strengthened when  
\(S\) or \(\Sigma\) satisfies a restriction theorem  
and to provide an estimate for the constant \(C(S, \Sigma)\).  
This result serves as a quantitative counterpart to a recent finding by the first and last author \cite{IM24}.  
In the setting of finite groups, the results also extend those of Matolcsi-Sz\"ucs and Donoho-Stark.  
\end{abstract}  



\maketitle

\tableofcontents

\section{Introduction} 

In this paper we establish a version of the uncertainty principle that roughly states that a function and its Fourier transform can not be arbitrarily
concentrated in a pair of sets of which one
has sufficiently small finite measure and the other
one is a sufficiently small neighborhood of
a set that satisfies the Fourier restriction property.

The Uncertainty Principle is one of the fundamental features of Fourier analysis. As a meta-theorem, it states that a function \( f \) and its Fourier transform \( \widehat{f} \) cannot be arbitrarily concentrated over a time and frequency domain simultaneously.
There are numerous mathematical formulations of the Uncertainty Principle and we refer to the book \cite{HJ1994} or the surveys \cite{FS1997,BD2006} for some of those formulations.
In this paper, we will deal with the Uncertainty Principle (\textit{UP}) on
locally compact abelian (\textit{LGA}) groups.  
The formulation we deal with is known as the {\em Annihilating Pair} property as well as the Amrein-Berthier-Benedicks UP in the continuous setting and Donoho-Stark UP in the discrete one. It states that a function and its Fourier transform can not be  supported in certain pairs of sets.

\smallskip

Let us now be more precise. First, in this paper, we adopt the general setting of \text{LGA} groups. Throughout, $G$ will be a \text{LGA} group with the Haar measure $\mu$ and $\widehat{G}$ its dual group, {\it i.e.}
the set of characters of $G$. Then $\widehat{G}$, the set of characters on $G$, which
is also a LCA group and we denote by $\widehat{\mu}$ its Haar measure.
Depending on the circumstances, it is more convenient to consider elements of $\widehat{G}$
as characters (continuous group   homomorphisms $\chi: G\to\{z\in\bC\ :\ |z|=1\}$) or as a set that parametrizes those functions.
The Fourier transform of $f\in L^1(G)$ is the function on $\widehat{G}$ defined by
$$
\widehat{f}(\chi)=\cF_G f(\chi)=c_G\int_G f(x)\overline{\chi(x)}\,\mbox{d}\mu(x),
\qquad\chi\in\hat G
$$
where $c_G$ is a normalization constant depending on $G$. The Fourier transform is an injection on $L^1(G)$. 

We will always assume that $c_G$ is chosen so that the 
Fourier transform extends to a unitary transform, i.e., $\|\widehat{f}\|_{L^2(\widehat{G})}=\|f\|_{L^2(G)}$.
Recall that if further 
$\widehat{f}\in L^1(\widehat{G})$, then $f=\cF^{-1}_G\widehat{f}$
with the inverse Fourier transform $\cF^{-1}_G$ defined by
$$
\cF^{-1}_Gg(x)=c_{\widehat{G}}\int_{\widehat{G}}g(\chi)\chi(x)\,\mbox{d}\widehat{\mu}(\chi),
$$ 
where $g=\hat f$. 
Note that $\cF^{-1}_G g=\overline{\cF_{\widehat G} \bar g}$ 
so that $\cF^{-1}_G$ also extends to a unitary transform on $L^2(\hat G)$
and that $f=\cF^{-1}_G\cF_Gf$ for $f\in L^2(G)$. See, for example, \cite{Rudin62,Terras99}. 

\vskip.125in 

To further fix normalizations, let us detail the two cases on which we will focus later on. 
\begin{itemize}[left=0pt]
\item[--] The {\em Discrete Fourier Transform}: Let $G=\widehat{G}=(\bZ/N\bZ)^d$ and $\mu=\widehat{\mu}$ is the counting measure. Here we identify the character $m\in (\bZ/N\bZ)^d$ with the function $e_m$ on $(\bZ/N\bZ)^d$ given by $e_m(n)=e^{2i\pi \scal{m,n}/N}$ so that
for $f: (\bZ/N\bZ)^d\to \mathbb{C}$
$$ \widehat{f}(m):=N^{-\frac{d}{2}}\sum_{n\in(\bZ/N\bZ)^d}f(n)e^{-2i\pi \scal{m,n}/N}. $$

The inverse Fourier transform formula is then given by 
$$
f(n):=N^{-\frac{d}{2}}\sum_{m\in(\bZ/N\bZ)^d}\widehat{f}(m)e^{2i\pi \scal{m,n}/N}, 
$$
and we have the Plancherel identity 
$$
\|f\|_{L^2(G)}^2:=\sum_{n\in (\bZ/N\bZ)^d}|f(n)|^2=\sum_{m \in (\bZ/N\bZ)^d} {|\widehat{f}(m)|}^2.
$$
Note that in this case, $c_G=c_{\widehat{G}}=N^{-\frac{d}{2}}$. 

\vskip.125in 

\item[--] The  {\em Continuous Fourier Transform:} Let $G=\bR^d$ equipped with the  Lebesgue
measure $dx$. It follows that $\widehat{G}=\bR^d$ and the Haar measure on $\widehat{G}$ is
  the Lebesgue measure $d\xi$. We identify $\xi\in \hat G=\bR^d$ with the character $e_\xi$ on $\bR^d$ given by $e_\xi(x)=e^{2i\pi \scal{x,\xi}}$. Then, for $f\in L^1(G)$, we define
$$
\widehat{f}(\xi):=\int_{\bR^d}f(x)e^{-2i\pi \scal{x,\xi}}\,\mbox{d}x,
$$
and if $\widehat{f}\in L^1(\bR^d)$, we have 
$$
f(x)=\int_{\bR^d}\widehat{f}(\xi)e^{2i\pi \scal{x,\xi}}\,\mbox{d}\xi,
$$
and then the Fourier transform extends from $L^1(\bR^d)\cap L^2(\bR^d)$ to $L^2(\bR^d)$
such that 
$$
\|f\|_{L^2(\bR^d)}^2:=\int_{\bR^d}|f(x)|^2\,\mbox{d}x
=\int_{\bR^d}|\widehat{f}(\xi)|^2\,\mbox{d}\xi:=\|\widehat{f}\|_{L^2(\bR^d)}^2.
$$
Here $c_G=c_{\widehat{G}}=1$.
\end{itemize}

\vskip.125in 

The main thrust of this paper is the following variant of Fourier uncertainty principle,  where the concentration of a function over a domain in $G$  is measured in terms of the size of the support of the function and its Fourier transform.

\begin{definition} Let $G$ be a locally compact abelian group quipped with a Haar measure $\mu$. Let $\hat G$ and $\hat\mu$
 be as described above. 
Let \( S \subset G \), \( \Sigma \subset \widehat{G} \), \( S^c = G \setminus S \), and \( \Sigma^c = \widehat{G} \setminus \Sigma \). Then 

\begin{itemize}
    \item the pair
\((S, \Sigma)\) is said to be a \emph{weak annihilating pair} if, for any \( a \in L^2(G) \) with \( \mathrm{supp}\, a \subset S \) and \( \mathrm{supp}\, \widehat{a} \subset \Sigma \), it follows that \( a = 0 \);

\item
  The pair $(S,\Sigma)$ is called   a \emph{strong annihilating pair} if there exists a constant $C(S,\Sigma)$ such that for every $a\in L^2(G)$
\begin{equation}
\label{eq:strongan}
\norm{a}_{L^2(G)}\leq C(S,\Sigma)\bigl(\norm{a}_{L^2(S^c)}
+\|\widehat{a}\|_{ L^2(\Sigma^c)}\bigr).
\end{equation}
\end{itemize}
The constant $C(S,\Sigma)$ is called the {\em annihilation constant} of $(S,\Sigma)$. 
\end{definition}

We will list several examples of such pairs below.
Note that a strong annihilating pair is also a weak one. The converse need not be true,
except if $L^2(G)$ is finite-dimensional, which is the case when $G$ is a finite abelian group. See, for example, \cite[chapter 3]{HJ1994}) (or the Appendix \ref{Appendix}) for the basic facts on this notion.

Before explaining why this notion has received a lot of attention in the last two decades, let us stress that annihilating pairs are strongly linked to the properties of the orthogonal projections $P_Sf=\mathbf{1}_Sf,Q_\Sigma=\cF^{-1}\mathbf{1}_\Sigma\cF$. More precisely, $(S,\Sigma)$ is strongly annihilating if and only if $\norm{P_SQ_\Sigma}_{L^2(G)\to L^2(G)}<1$ and then 
$$
C(S,\Sigma)=\dfrac{1}{\sqrt{1-\norm{P_SQ_\Sigma}_{L^2(G)\to L^2(G)}^2}}.
$$ 
This explains why most results in this area rely on Hilbert space methods. Even though some results
extend to $L^p$ space, they mainly rely
on $L^p\to L^p$ properties of the map $P_SQ_\sigma$. 
The goal of this paper is to show that using a larger scale of $L^p$-spaces and $L^p\to L^q$-properties of the Fourier transform, such as the Hausdorff-Young inequality and Fourier restriction estimates, we can obtain new insights on annihilating pairs, even in the $L^2$ setting.
We will detail this further in the next section.

\vskip.125in 

Let us now explain how annihilating pairs appear in many aspects of applied sciences. 
To start, let us mention that in practice, the support may not be the right paradigm as a signal $s$ will not be detected below some threshold or if it does not carry enough energy. A way to formulate this mathematically is to introduce the following notion. 

\begin{definition}\label{time-frequency-concenteration} Let $G$ be a locally compact abelian group, $S\subset G$, $\Sigma\subset\widehat{G}$ and $\eps_T,\eps_\Omega>0$. We will say that a signal $f\in L^2(G)$ is {\em $\eps_T$-timelimited} to $S$ if 
$$
\|f\|_{L^2(G\setminus S)}\leq\eps_T\|f\|_{L^2(G)},
$$
and {\em $\eps_\Omega$-bandlimited} to $\Sigma$ if 
$$
\|\widehat{f}\|_{L^2(\widehat{G}\setminus \Sigma)}\leq\eps_\Omega\|\widehat{f}\|_{L^2(\widehat{G})}.
$$    
\end{definition}

This notion has been introduced by Landau, Pollak and Slepian in the case $G=\R$, $S=S_T=[-T,T]$ and $\Sigma=\Sigma_\Omega=[-\Omega,\Omega]$. This idea is closely linked to annihilating pairs as follows: If $(S_T,\Sigma_\Omega)$ is a strong annihilating pair, then any $f\not=0$ that is $\eps_T$-timelimited to $S$ and $\eps_\Omega$-bandlimited to $\Sigma$ will satisfy
$$
\|f\|_{L^2(\R)}\leq C(S_T,\Sigma_\Omega)\bigl(\|f\|_{L^2(\R\setminus S_T)}+\|\widehat{f}\|_{L^2(\R\setminus \Sigma_\Omega)}\bigr)
\leq  C(S_T,\Sigma_\Omega)(\eps_T+\eps_\Omega)\|f\|_{L^2(\R)}.
$$
This shows that $\eps_T,\eps_\Omega>0$ can not be taken arbitrarily small as one requires
\begin{equation}
    \label{eq:concentration}
\eps_T+\eps_\Omega\geq \frac{1}{C(S_T,\Sigma_\Omega)}.
\end{equation}

In the opposite direction, if $\epsilon_T, \epsilon_\Omega > 0$ are small but sufficiently large, the seminal work of Landau, Pollak, and Slepian proves that prolate spheroidal wave functions satisfy these limiting properties, provided that $\Omega T$ is not too small with respect to the Nequist rate.
 The key here is the 
investigation of the eigenvalue behavior of   the compact self-adjoint and positive definite operators $P_SQ_\Sigma P_S$, {\it time-frequency limiting operators}. Let us mention that its eigenvalues also play an important role in random matrix theory.

\medskip 


\smallskip

Another application coming from signal processing is the following. Let $G$ be a finite group, say $G=(\bZ/N\bZ)^d$ with $d=1$ or $2$. We call an element $s$ of $L^2(G)$ a signal to make matters more explicit. We will say that a signal $s$ is $t$-sparse if its  spectrum 
has only $t$-elements, i.e the support of the Fourier transform of $s$ has size $t$, and in practice, many natural signals can reasonably considered to be $t$-sparse for some $t\ll |G|=N^d$
(or at least $\eps_\Omega$-concentrated to a set of size
$t\ll N^d$). Unfortunately, not all components of the signal  $s$
  can generally be measured, or they may be so heavily corrupted by noise that is better to discard them. This phenomenon occurs in many real-world signals (such as images, audio, and medical signals), which often has sparsity when represented in an appropriate basis (e.g., wavelets, Gabor bases, or Fourier transforms).
We are thus given a set $T\subset G$ of known locations on which the measurement of $s$
is sufficiently reliable, that is, our measurement
is $s$ restricted to $T$. We thus want to recover $s$ from its restriction to
$T$ and the knowledge that it is $t$ sparse.
The first question we thus ask is whether such an $s$
may be uniquely determined. Assume that $s_1,s_2$
are both $t$-sparse and agree on $T$. 
If we make the stronger assumption that we actually know
exactly the support $\Upsilon$ of $\widehat{s_1},\widehat{s_2}$, and if
$S=G\setminus T$ and $\Sigma=\widehat{G}\setminus\Upsilon$
form an annihilating pair 
then, as $s=s_1-s_2$ vanishes on $S$ and its Fourier transform
vanishes on $\Sigma$ (since both $s_1$ and $s_2$ do)
then $s=0$, that is $s_1=s_2$. 
Now if we only know the size $t$ of the support $\widehat{s_1},\widehat{s_2}$, then $\widehat{s}=
\widehat{s_1}-\widehat{s_2}$ has support of size at most $2t$.
To obtain uniqueness, we are now asking whether, given $S=G\setminus T$
and any set $\Sigma$ of size $|\widehat G|-2t$, $(S,\Sigma)$ forms an annihilating pair. 

In this case, the only information we have is on the size of $\Sigma$. The first question we ask is
whether $S$ can be large or small.
This is a well-known result by   Matolcsi and Szucks \cite{MS73}:
if $S,\Sigma\subset (\bZ/N\bZ)^d$ and $|S||\Sigma|<N^d$, 
then $(S,\Sigma)$ is a (strong) annihilating pair.
The result is optimal in the sense that, without further information on $S$ than its size, there exists $S,\Sigma$
such that $|S||\Sigma|=N^d$ that are not annihilating pairs.
This result was rediscovered and further developed by Donoho and Stark \cite{DS89}. The constant $C(S,\Sigma)$
has been computed for any finite abelian group by Ghober and the second author \cite{GJ2011}. This result for $G=(\bZ/N\bZ)^d$ is stated as follows: 
$$ \norm{f}_{\ell^2(\bZ/N\bZ)^d}\leq 
\left(1+\frac{1}{1-\sqrt{\frac{|S||\Sigma|}{N^d}}}\right)
\bigl(\norm{f}_{\ell^2((\bZ/N\bZ)^d\setminus S)}+
\|\widehat{f}\|_{\ell^2((\bZ/N\bZ)^d\setminus\Sigma)}\bigr).
$$
Some improvements are possible in specific cases.
For instance, when $G=\bZ/p\bZ$, $p$ a prime number, then 
the optimal condition is $|S|+|\Sigma|<p+1$
as proved by Tao \cite{TaoUncertainty}.
Further results of that flavor can be found in
\cite{BG2013,Meshulam06}.  

This then leads to the question of finding an efficient algorithm for reconstruction
{\it i.e.} to find a signal $s\in L^2(G)$
such that $s$ agrees with the measurement $\mu$
on $T$, $s(j)=\mu(j)$ when $j\in T$,  and such that
$\widehat{s}$ is supported in $\Upsilon$.
If $\Upsilon$
were given, this could be done in many ways, including
with various versions of Kaczmarz's algorithm \cite{Kaczmarz}. This amounts to iterating the following operation, also known as Iterative Projection Method. 
\begin{itemize}
    \item {\em Initialise} $s_0=\mu\mathbbm{1}_T+0\mathbbm{1}_{G\setminus T}$,
    \item Repeat
    \begin{itemize}
        \item Compute $\cF_G[s_j]$ and keep
        only coefficients in $\Upsilon$
        and invert the Fourier transform: 
        $\nu=\cF_G^{-1}[\cF_G[s_j]\mathbf{1}_{\Upsilon}]$
        \item Update $s_j$ by keeping the measurements on $T$ and replace with the coefficients of $\nu$ elsewhere
        $$
        s_{j+1}=\mu\mathbbm{1}_T+\nu\mathbbm{1}_{G\setminus T}.
        $$
    \end{itemize}
    \item Stop once the iteration $\nu$ is small enough or after a prefixed maximum number of iterations.
\end{itemize}

This algorithm is easy to programm and usually runs fast and provides us with the solution of the following program:
$$
\mathrm{argmin}\{\|s-\mu\|_{\ell^2(T)}\,:\supp\widehat{s}\subset\Upsilon\}
$$
The drawback is that it requires the knowledge of the set $\Upsilon$.
If we only know that the size of $\Upsilon$ is at most $t$,
our program is slightly different as we want to solve
$$
\mathrm{argmin}\{\|\widehat{s}\|_0\,: s=\mu\mbox{ on }T\}
$$
where $\|\widehat{s}\|_0=|\supp\widehat{s}|$.
The naive 
idea would be to consider any set $\Upsilon$ 
of size at most $t$
and run the previous algorithm and then repeat this until
one finds the appropriate $\Upsilon$. 
This may unfortunately
imply the run over all possible sets which is $NP$-complete and is thus unpracticable.
Fortunately, an alternative is sometimes possible.
It turns out that when $S=G\setminus T$ and 
$\Sigma=\widehat{G}\setminus\Upsilon$
form a strong annihilating pair and further that
the matrix $I-P_SQ_\Sigma$ is not only invertible,
but also sufficiently
well-conditioned, then a convex relaxation 
$$
\mathrm{argmin}\{\|\widehat{s}\|_1\,: s=\mu\mbox{ on }T\}
$$
will produce the solution.
This is at the heart of {\em compressed sensing}
which experienced a burst of popularity a decade ago, {\it see e.g.}
\cite{CR05,CRT06} and the books \cite{AH21,FR2013} for more on the subject.

\smallskip

Let us now move to the continuous setting. 
A first observation is that many constant coefficient
PDEs may be solved via the Fourier transform.
For instance, if $v$ is a solution of the
free heat equation
$$
\left\{\begin{matrix} \partial_tv+\frac{1}{(2\pi)^2}\Delta_xv=0\\
v(x,0)=v_0(x) \end{matrix}\right.
$$
with $v_0\in L^2(\R^d)$, then 
for any $t\in \Bbb R^{+}$ 
$$
v(x,t)=\int_{\R^d}e^{-\pi t|\xi|^2}\widehat{v_0}(\xi)
e^{2i\pi\scal{x,\xi}}\,\mbox{d}\xi=
\cF^{-1}[e^{-\pi t|\cdot|^2}\widehat{v_0}(\cdot)](x).
$$
In particular, if $\Sigma=B(0,a)$ is  a ball of radius $a$ in $\Bbb R^d$,  and $S\subset \Bbb R^d$,    such that
$(S,\Sigma)$ is a strong annihilating pair with constant
$C(S,\Sigma)$, then
for any fixed $t\in \Bbb R^+$
 
$$
\norm{v(x,t)}_{L^2(\R^d\setminus S)}\geq C(S,\Sigma)
e^{-\pi a^2t}\norm{v_0}_{L^2(\R^d)}.
$$
 
This can be seen as a quantitative unique continuation property of the heat equation.

Moreover, in this case, the sets $S$
for which $(S,\Sigma)$ is a strong annihilating pair
have been characterized by Logvinenko-Sereda \cite{LS74}
({\it see} \cite{HJ1994} for further references on this work). They are the so-called dense sets, meaning that there exists $0<\gamma<1$ and $r>0$
such that for every $x\in\R^d$, $|S\cap B(x,r)|\geq\gamma|B(x,r)|$.
Two decades ago, Kovrijkine \cite{K01} proved the estimate $C(S,\Sigma)\leq\left(\dfrac{C}{\gamma}\right)^{cab+1}$ which is optimal concerning the behavior with respect to the parameters $a,b$ and $\gamma$
(up to the explicit numerical constants $c,C$). Recently, this result has found applications in control theory
as Egidi-Veselic \cite{EV18} and Wang, Wang, Zhang, Zhang \cite{WWZZ19} independently showed that this result implies that the heat equation on the full space $\R^d$ is null controllable from a set $\Omega$ if and only if $\Omega$ is relatively dense.

A second important result is that if $S,\Sigma$
are sets of finite measure in $\R^d$, then
Benedicks \cite{Be85} showed that $(S,\Sigma)$
is a weak annihilating pair while Amrein-Berthier \cite{AB77} showed that it is actually a strong
annihilating pair (see \cite{BD2006} for an argument showing how this is implied by Benedicks's result).
Their original motivation was the analysis of the joint measurement of incompatible observables in quantum mechanics. The constant $C(S,\Sigma)$ was later shown by Nazarov \cite{N93} to be of the form $Ce^{C|S||\Sigma|}$
when $G=\R$. The result was extended to $G=\R^d$, $d\geq 2$
by the second author \cite{Ja07}, though in this
case it is conjectured that the right behavior should be
$Ce^{C(|S||\Sigma|)^{1/d}}$. In \cite{Ja07}, the conjecture was also verified to hold when one of $S$ or $\Sigma$ is convex.

There are many other examples. Let us mention that
strong annihilating pairs play a key role in the investigation of Anderson localization by Shubin, Vakilian, Wolff \cite{SVW98}. A further important
example is Bourgain-Dyatlov's Fractal Uncertainty Principle \cite{BoDy2018}, which
has been successfully applied to problems in quantum chaos; {\it see} the survey \cite{Dyatlov-FUP} and references therein.

\medskip

The remaining of this paper is organized as follows: In the next section, we present the necessary machinery and give a general statement on locally compact abelian groups. We then specify this to two settings. In Section 3, we will deal with the case of finite abelian groups. We conclude in Section 4 with the case of $G=\R^d$.

\section{The main result in an abstract setting}

\subsection{Restriction Estimates}

We will write $\cS(G)\subset L^1(G)$ a set of functions that is dense in every $L^p(G)$ space, $1\leq p<+\infty$
and such that, for $f\in\cS(G)$, $\widehat{f}$ is continuous.
 
\begin{definition} Let $G$ be a LCA group, $\widehat{G}$ its dual group, and 
$1 \leq p \leq q \leq \infty$.
A set $\Sigma\subset \widehat{G}$ with  $\widehat{m}(\Sigma)>0$ is said to satisfy a {\em $(p,q)$-restriction estimate} with constant
$\rho_{p,q}(\Sigma)$ if, for all $f\in\cS(G)$,
\begin{equation}
    \label{eq:resleb}
\left(\int_\Sigma|\widehat{f}(\xi)|^q\,\mathrm{d}\widehat{m}(\xi)\right)^{\frac{1}{q}}
\leq \rho_{p,q}(\Sigma)\left(\int_G|f(x)|^p\,\mathrm{d}m(x)\right)^{\frac{1}{p}}. 
\end{equation}
When $p$ or $q=+\infty$, we replace the norms with the supremum norms.
\end{definition}

The equation \eqref{eq:resleb} implies that the Fourier transform extends into a continuous operator $L^p(G)\to L^q(\Sigma)$.

The point of restriction estimates is that, under suitable geometric conditions, the restriction estimate
may hold for a larger set of $(p,q)$'s or with better constants.
 In the estimation, 
ideally, one wishes for the constant $\rho_{p,q}(\Sigma)$ to be small.
 
\smallskip

For instance, in the Euclidean case $G=\R^d$, the restriction estimate is usually stated for a smooth
hypersurfaces $\mathbb{S}$
endowed with its surface measure $\nu$. The $L^q$-norm in \eqref{eq:resleb} is then replaced by the $L^q(\mathbb{S},\nu)$-norm.
However, if $\Sigma=\mathbb{S}_\delta=
\{x:\exists y\in \mathbb{S}, |x-y|<\delta\}$
a $\delta$-neighborhood of $\mathbb{S}$,  then such an estimate is equivalent
to a restriction estimate in the sense of \eqref{eq:resleb} for all sufficiently small $\delta$,
provided the constant is of the form
$$
\rho_{p,q}=C_{p,q}(\mathbb{S})\delta^{1/q}.
$$ 
In other words,
\begin{equation}
    \label{eq:rescont}
\left(\int_{\mathbb{S}_\delta}|\widehat{f}(\xi)|^q\,\mathrm{d}\xi\right)^{\frac{1}{q}}
\leq C_{p,q}(\mathbb{S})\delta^{1/q}\left(\int_{\R^d}|f(x)|^p\,\mathrm{d}x\right)^{\frac{1}{p}}
\qquad \forall 0<\delta<\delta_0,\quad\forall f\in\cS(\R^d).
\end{equation}
We will give some examples in Section \ref{Sec:exdis}.  

\vskip.125in 

In the case when $G$ is a finite abelian group, every set $\Sigma$ satisfies $(p,q)$-restriction estimate \eqref{eq:resleb}, but the constant may be very large. The notion of restriction in the finite setting estimates, introduced in \cite{MT04}, is the following: We say that a subset $\Sigma\subset \widehat{G}$  satisfies the $(p,q)$-restriction estimate if there is a constant $C_{p,q}(\Sigma)$ such that for any function $f$ on $G$
\begin{equation} \label{restrictionequation} 
\left( \frac{1}{|\Sigma|} \sum_{\chi \in \Sigma} |\widehat{f}(\chi)|^q \right)^{\frac{1}{q}}
\leq \frac{C_{p,q}(\Sigma)}{|G|^{\frac{1}{2}}}
\left(\sum_{x \in G} |f(x)|^p \right)^{\frac{1}{p}}. 
\end{equation} 
This is of course the same as \eqref{eq:resleb} with
\begin{equation}
\label{eq:goodrestdis}
\rho_{p,q}(\Sigma)=C_{p,q}(\Sigma)\frac{|\Sigma|^{\frac{1}{q}}}{|G|^{\frac{1}{2}}}.
\end{equation}
For certain families $\mathcal{F}$ of $\Sigma$ and certain exponents $p, q$, $C_{p,q}$ is a constant that depends on the family $\mathcal{F}$ but not on the particular element $\Sigma$ of the family.
 
There is a vast literature on the subject; see,  e.g., \cite{MT04,HW18,IK08,IK10,IK10b,IKL17}. For further detail, see Section \ref{Sec:exdis}.

\subsection{The main result}
In the following, we assume that $G$ is a LCA group equipped with a Haar measure $m$ and $\widehat{G}$ is its dual group with Haar measure $\hat m$. 
\begin{theorem}\label{th:IMquant}
  Let $1\leq p\leq 2\leq q$. Let $S\subset G$ and $\Sigma\subset\widehat{G}$ be two sets of finite measure.
  Assume further that $\Sigma$ satisfies a $(p,q)$-restriction estimate \eqref{eq:resleb} with constant $\rho_{p,q}(\Sigma)$.
  
  Assume that $m(S)$ and $\widehat{m}(\Sigma)$ are small enough to satisfy
  \begin{equation}
      \label{eq:condthImquant}
  \rho_{p,q}(\Sigma)m(S)^{\frac{1}{p}-\frac{1}{2}}\widehat{m}(\Sigma)^{\frac{1}{2}-\frac{1}{q}} < 1.
 \end{equation}
  Then $(S,\Sigma)$ is a strong annihilating pair; i.e. 
$$
\|f\|_{L^2(G)}\leq A_{ann}(S,\Sigma)\bigl(\|f\|_{L^2(G\setminus S)}+\|\widehat{f}\|_{L^2(\widehat{G}\setminus \Sigma)}\bigr).
$$
with
$$
A_{ann}(S,\Sigma)=\frac{1}{1-\rho_{p,q}(\Sigma)m(S)^{\frac{1}{p}-\frac{1}{2}}\widehat{m}(\Sigma)^{\frac{1}{2}-\frac{1}{q}}}.
$$
\end{theorem}
\begin{remark}
    Any set $\Sigma$ of finite measure satisfies the $(1,\infty)$ restriction estimate with constant 
    $\rho_{1,\infty}(\Sigma)=c_G$, a constant independent of $\Sigma$.
    When $G$ is a finite group, we recover Matolcsi-Suck's Uncertainty Principle in a quantitative form, any pair of sets $(S,\Sigma)$ such that $|S|\,|\Sigma|<|G|$ is a strong annihilating pair,
    and we also recover the constant from \cite{GJ2011}.

    Note also that the statement is void for $p=q=2$, though any set $\Sigma$ of finite measure also
    satisfies the $(2,2)$-restriction estimate with constant 
    $\rho_{2,2}(\Sigma)=1$, again independent of $\Sigma$.

    Note also that when $G=(\bZ/N\bZ)^d$,  the fact that
    $(S,\Sigma)$ is a weak (thus strong) annihilating pair
    was previously proven by the first and last author. It is then better to write $\rho_{p,q}(\Sigma)$
    in the form given by \eqref{eq:goodrestdis} and the condition is then
    $$
   |\Sigma||S|^{\frac{2}{p}-1}<\frac{|G|}{C_{p,q}(\Sigma)^2}.
   $$
   The improvement is that the exponent of $|S|$ is $<1$ so that, in the presence of a good bound on
   $C_{p,q}(\Sigma)^2$, this condition is much less restrictive than in Matolcsi and Sucks when $|G|$ is large.
   Also, the novelty here is that we compute the annihilation constant
    $$
    A_{ann}(S,\Sigma)=1+\frac{1}{1-C_{p,q}(\Sigma)\sqrt{\frac{|S|^{\frac{2}{p}-1}|\Sigma|}{|G|}}}.
    $$
\end{remark}

\begin{proof}[Proof of Theorem 
\ref{th:IMquant}]
To simplify notation, we write $\rho_{p,q}$ instead of $\rho_{p,q}(\Sigma)$. H\"older's inequality implies that
$$
\|\widehat{\mathbf{1}_S f}\|_{L^2(\Sigma)} \leq \widehat{m}(\Sigma)^{\frac{1}{2}-\frac{1}{q}}\|\widehat{\mathbf{1}_S f}\|_{L^q(\Sigma)}
  \leq  \rho_{p,q}(\Sigma) \widehat{m}(\Sigma)^{\frac{1}{2}-\frac{1}{q}}\|f\|_{L^p(S)}
  $$
where we used the restriction estimate on $\Sigma$ to obtain the rightmost inequality. Applying the H\"older inequality for $\|f\|_{L^p(S)}$ one more time   gives
$$
\|\widehat{\mathbf{1}_S f}\|_{L^2(\Sigma)} 
\leq  \rho_{p,q}(\Sigma)\widehat{m}(\Sigma)^{\frac{1}{2}-\frac{1}{q}}m(S)^{\frac{1}{p}-\frac{1}{2}}\|f\|_{L^2(S)}.
  $$
To simplify notation till the end of the proof, write
this in the form $\|\widehat{\mathbf{1}_S f}\|_{L^2(\Sigma)} \leq A\|f\|_{L^2(S)}$, where $A=\rho_{p,q}(\Sigma)\widehat{m}(\Sigma)^{\frac{1}{2}-\frac{1}{q}}m(S)^{\frac{1}{p}-\frac{1}{2}}$,  and note that our hypothesis is that $A<1$.
The remaining of the proof is standard machinery. First
$$
\|\widehat{\mathbf{1}_S f}\|_{L^2(\widehat{G}\setminus\Sigma)}\geq  
\|\widehat{\mathbf{1}_S f}\|_{L^2(\widehat{G})}
-\|\widehat{\mathbf{1}_S f}\|_{L^2(\Sigma)} \geq 
(1-A)\|\mathbf{1}_S f\|_{L^2(G)}
$$
with Parseval.  
This implies that 
$$\|f\|_{L^2(S)}\leq  \frac{1}{1-A}\|\widehat{\mathbf{1}_S f}\|_{L^2(\widehat{G}\setminus\Sigma)}
$$


Note that  we can write 
$$
\|\widehat{\mathbf{1}_S f}\|_{L^2(\widehat{G}\setminus\Sigma)}\leq
\|\widehat{f}\|_{L^2(\widehat{G}\setminus\Sigma)}
+\|\widehat{\mathbf{1}_{G\setminus S} f}\|_{L^2(\widehat{G}\setminus\Sigma)}
$$
and the last quantity is bounded by
$$
\|\widehat{\mathbf{1}_{G\setminus S} f}\|_{L^2(\widehat{G}\setminus\Sigma)}\leq
\|\widehat{\mathbf{1}_{G\setminus S} f}\|_{L^2(\widehat{G})}=\|f\|_{L^2(G\setminus S)}
$$
using Parseval. All together, by writing  $f= f{\bf 1}_S + f {\bf 1}_{G\setminus S}$ and the fact that $0< A<1$, we obtain that
$$
\|f\|_{L^2(G)}  \leq  \|f\|_{L^2(S)}+\|f\|_{L^2(G\setminus S)}
\leq \frac{1}{1-A}
\left(
\|\widehat{f}\|_{L^2(\widehat{G}\setminus\Sigma)}+
\|f\|_{L^2(G\setminus S)}\right).
$$
The proof follows by replacing $A$ by its value and  $A_{ann}(S,\Sigma) = \cfrac{1}{1-A}$.
\end{proof}

\section{Examples when $G$ is a finite Abelian group}
\label{Sec:exdis}

In this section $G$ will be a finite Abelian group.

\subsection{A second result}

The first observation is that, when $G$ is finite,
we are not limited to $q\geq 2$ for the use
of a restriction estimate.

Recall that when $G$ is a finite abelian group, then so is $\widehat{G}$
and that Haar measure on those groups is the counting measure. We will thus write
$|S|,|\Sigma|$ for the cardinality/measure of sets $S\subset G$, $\Sigma\subset\widehat{G}$.
Recall also that $|\widehat{G}|=|G|$.
Finally, the Fourier transform normalization constants are $c_G=c_{\widehat{G}}=|G|^{-1/2}$.

\begin{theorem}\label{th:mainThm2}
    Let $G$ be a finite abelian group.
    Let $1\leq p,q\leq 2$. Let $\Sigma\subset\widehat{G}$ that satisfies the $(p,q)$-restriction estimate \eqref{restrictionequation} with the constant $\rho_{p,q}(\Sigma)$ given as in \eqref{eq:goodrestdis}. 
Let $S\subset G$ be such that 
\begin{equation}
    \label{eq:condonr}
    \rho_{p,q}(\Sigma)|S|^{\frac{1}{p}}<|G|^{\frac{1}{q}-\frac{1}{2}}
\end{equation}
    Then $(S,\Sigma)$ is a strong annihilating pair with constant 
      $$
A_{ann}(S,\Sigma)=
1+\frac{|S|^{\frac{1}{2}} ~ |\widehat{G}\setminus\Sigma|^{\frac{1}{q}-\frac{1}{2}}} {|G|^{\frac{1}{q}-\frac{1}{2}}-\rho_{p,q}(\Sigma)|S|^{\frac{1}{p}}}.
$$ That is, 
    for every $f\,:G\to\bC$,
    \begin{equation}
    \|f\|_{L^2(G)}\leq A_{ann}(S,\Sigma)\bigl(\|f\|_{L^2(G\setminus S)}+\|\widehat{f}\|_{L^2(\widehat{G}\setminus \Sigma)}\bigr). 
\label{eq:thAIfinite}
    \end{equation}
 
\end{theorem} 
\begin{remark}
Notice that a simple application of Hölder's inequality shows that if $\Sigma$ satisfies the  
$(p,q)$-restriction estimate with constant $\rho_{p,q}(\Sigma)$, then  
it also satisfies the $(r,s)$-restriction estimate for any $r \geq p$ and $s \leq q$ with constant $\rho_{r,s}(\Sigma)$, where  
\begin{equation}
    \label{eq:changerestriction}
\rho_{r,s}(\Sigma)=|\Sigma|^{\frac{1}{s}-\frac{1}{q}}|G|^{\frac{1}{p}-\frac{1}{r}}\rho_{p,q}(\Sigma).
\end{equation}
In order to compare this result with Theorem \ref{th:IMquant}, we write $\rho_{p,q}(\Sigma)$
in the form given by \eqref{eq:goodrestdis} {\it i.e.} 
$\rho_{p,q}(\Sigma)=C_{p,q}(\Sigma)\dfrac{|\Sigma|^{\frac{1}{q}}}{|G|^{\frac{1}{2}}}$.
Then the condition in Theorem \ref{th:IMquant} is, for $q\geq 2$,
$$
\rho_{p,q}(\Sigma)|S|^{\frac{1}{p}-\frac{1}{2}}|G|^{\frac{1}{2}-\frac{1}{q}}<1
\Longleftrightarrow
C_{p,q}(\Sigma)|S|^{\frac{1}{p}-\frac{1}{2}}|\Sigma|^{\frac{1}{q}}<|G|^{\frac{1}{q}}.
$$
If we apply the $(r,s)$-restriction estimate instead, and use \eqref{eq:changerestriction}, the condition becomes
$$
\left(\frac{|G|}{|S|}\right)^{\frac{1}{p}-\frac{1}{r}}\left(\frac{|\Sigma|^2}{|G|}\right)^{\frac{1}{s}-\frac{1}{q}}
C_{p,q}(\Sigma)|S|^{\frac{1}{p}-\frac{1}{2}}|\Sigma|^{\frac{1}{q}}
<|G|^{\frac{1}{q}}.
$$
This shows that increasing $r$ will add a factor $>1$ and lead to a more restrictive condition on $S$.
The situation is more complicated for $s$. If $\Sigma$ is ``big", $|\Sigma|\gg\sqrt{|G|}$, the condition
becomes more restrictive, while for $|\Sigma|\ll\sqrt{|G|}$, it is best to replace $q$ by the smallest $s$
possible, that is $s=2$, leading to the condition
$$
|\Sigma|\leq|G|^{\frac{1}{2}}\quad\mbox{and}\quad C_{p,q}(\Sigma)|S|^{\frac{1}{p}-\frac{1}{2}}|\Sigma|^{1-\frac{1}{q}}
<|G|^{\frac{1}{2}}.
$$

When $q=2$, there is no improvement possible in Theorem \ref{th:IMquant} and the condition is
$$
C_{p,2}(\Sigma)|S|^{\frac{1}{p}-\frac{1}{2}}|\Sigma|^{\frac{1}{2}}
<|G|^{\frac{1}{2}}\Longleftrightarrow
C_{p,2}(\Sigma)^2|S|^{\frac{2}{p}-1}|\Sigma|<|G|.
$$
However, when $q\leq 2$, Theorem \ref{th:mainThm2} applies and we obtain 
$$
C_{p,q}(\Sigma)^q|S|^{\frac{q}{p}}|\Sigma| <|G|.
$$
When $q=2$, this condition is 
$$
C_{p,2}(\Sigma)^2|S|^{\frac{2}{p}}|\Sigma|<|G|
$$
which is more restrictive than the one obtained from Theorem \ref{th:IMquant}.
However, if the $(p,2)$-restriction estimate holds, so does the $(p,q)$ one for $q\leq 2$
with constant $\rho_{p,q}(\Sigma)=|\Sigma|^{\frac{1}{q}-\frac{1}{2}}\rho_{2,q}$.
We can then apply Theorem \ref{th:mainThm2} as soon as
$$
C_{p,2}(\Sigma)^q|S|^{\frac{q}{p}}|\Sigma|^{2-\frac{q}{2}}
    <|G|.
$$
The condition is more restrictive for $\Sigma$ but becomes less restrictive for $S$
as soon as $\dfrac{q}{p}<\dfrac{2}{p}-1$, that is, if $q\leq 2-p$.
For instance, for $q=1$, this reads $C_{p,2}(\Sigma)|S|^{\frac{1}{p}}|\Sigma|^{\frac{3}{2}}<|G|$.
\end{remark}

\begin{proof}[Proof of Theorem \ref{th:mainThm2}]
As $\Sigma$ satisfies the restriction estimate, we have
\begin{equation}
    \label{eq:ph1}
    \|\widehat{\mathbf{1}_{S}f}\|_{L^q(\Sigma)}\leq \rho_{p,q}\|\mathbf{1}_{S}f\|_{L^p(G)}
\leq \rho_{p,q}|S|^{\frac{1}{p}-\frac{1}{2}}\|f\|_{L^2(S)}
\end{equation}
with H\"older's inequality.
Next, we use Hausdorff-Young's inequality for $\cF_G^{-1}$.
We write $\mathbf{1}_{S}f=\cF_G^{-1}\bigl[\cF_G[\mathbf{1}_{S}f]\bigr]$ thus
$$
\|\mathbf{1}_{S}f\|_{L^\infty(G)}\leq |G|^{\frac{1}{q}-\frac{1}{2}}\|\widehat{\mathbf{1}_{S}f}\|_{L^q(\widehat{G})}.
$$

\smallskip

Using H\"older's inequality, we obtain
\begin{eqnarray*}
\|f\|_{L^2(S)}=
 \|\mathbf{1}_{S}f\|_{L^2(G)}&\leq& |S|^{\frac{1}{2}}\|\mathbf{1}_{S}f\|_{L^\infty(G)}\\
 &\leq&|G|^{\frac{1}{2}-\frac{1}{q}}|S|^{\frac{1}{2}}\|\widehat{\mathbf{1}_{S}f}\|_{L^q(\widehat{G})}\\
&\leq&|G|^{\frac{1}{2}-\frac{1}{q}}|S|^{\frac{1}{2}}\bigl(\|\widehat{\mathbf{1}_{S}f}\|_{L^q(\Sigma)} 
+\|\widehat{\mathbf{1}_{S}f}\|_{L^q(\widehat{G}\setminus\Sigma)}\bigr).
\end{eqnarray*}

Combining this with \eqref{eq:ph1}, we get
$$
\frac{|G|^{\frac{1}{q}-\frac{1}{2}}}{|S|^{\frac{1}{2}}}\|f\|_{L^2(S)}
\leq \rho_{p,q}\frac{|S|^{\frac{1}{p}}}{|S|^{\frac{1}{2}}}\|f\|_{L^2(S)}
+\|\widehat{\mathbf{1}_{S}f}\|_{L^q(\widehat{G}\setminus\Sigma)}.
$$
As we assume that 
$$
\rho_{p,q}|S|^{\frac{1}{p}}<|G|^{\frac{1}{q}-\frac{1}{2}}
$$
then we reformulate this as
\begin{equation}
    \label{eq:ph2}
\|f\|_{L^2(S)}\leq\frac{|S|^{\frac{1}{2}}}{|G|^{\frac{1}{q}-\frac{1}{2}}-\rho_{p,q}|S|^{\frac{1}{p}}}\|\widehat{\mathbf{1}_{S}f}\|_{L^q(\widehat{G}\setminus\Sigma)}.
\end{equation}
Next, as in the proof of Theorem \ref{th:IMquant}, we unravel the norm of $\widehat{\mathbf{1}_{S}f}$ by writing
\begin{equation}
    \label{eq:ph3}
\|\widehat{\mathbf{1}_{S}f}\|_{L^q(\widehat{G}\setminus\Sigma)}
\leq |\widehat{G}\setminus\Sigma|^{\frac{1}{q}-\frac{1}{2}}\|\widehat{\mathbf{1}_{S}f}\|_{L^2(\widehat{G}\setminus\Sigma)}
\leq |\widehat{G}\setminus\Sigma|^{\frac{1}{q}-\frac{1}{2}}\bigl(\|\widehat{f}\|_{L^2(\widehat{G}\setminus\Sigma)}+\|f\|_{L^2(G\setminus S)}\bigr).
\end{equation}
Finally writing $\|f\|_{L^2(G)}\leq \|f\|_{L^2(G\setminus S)}+\|f\|_{L^2(S)}$ and bounding
the last term with \eqref{eq:ph2}-\eqref{eq:ph3} gives \eqref{eq:thAIfinite} with
$$
A_{ann}(S,\Sigma)=
1+\frac{|S|^{\frac{1}{2}}|\widehat{G}\setminus\Sigma|^{\frac{1}{q}-\frac{1}{2}}}{|G|^{\frac{1}{q}-\frac{1}{2}}-\rho_{p,q}|S|^{\frac{1}{p}}}.
$$
as claimed.
\end{proof}

\begin{remark}
    The reader may check that one could use Hausdorff-Young's inequality for an other exponent $r\geq 2$.
    The optimal constants for Hausdorff-Young's inequality on finite abelian groups are known as well as their optimizers, {\it see} \cite{GR2010}. Doing so, one would obtain a more restrictive condition on $S$ 
    $$
 \left(\frac{|G|}{|S|}\right)^{\frac{1}{r}}   |S|^{\frac{1}{p}}|\Sigma|^{\frac{1}{q}}
<\frac{|G|^{\frac{1}{2}+\frac{1}{q}}}{C_{p,q}}
$$
and a larger constant
$$
A_{ann}(S,\Sigma)=1+|G|^{\frac{1}{2}+\frac{1}{q}}\frac{\left(\frac{|S|}{|G|}\right)^{\frac{1}{2}-\frac{1}{r}}\left(1-\frac{|\Sigma|}{|G|}\right)^{\frac{1}{q}-\frac{1}{2}}}{|G|^{\frac{1}{2}+\frac{1}{q}}
-C_{p,q}\left(\frac{|G|}{|S|}\right)^{\frac{1}{r}}|S|^{\frac{1}{p}}|\Sigma|^{\frac{1}{q}}}.
$$
\end{remark}

The remainder of this section is devoted to some examples.

\subsection{$\Lambda_q$-sets}

In this section, $G$ will be either a finite or a compact group. 
Recall that elements of $\widehat{G}$ are continuous functions $G\to\{z\in\bC\,:\ |z|=1\}$.
Note that in both cases, if $q>2$, then $L^q(G)\subset L^2(G)$. In particular, $\widehat{G}\subset L^2(G)$.
Moreover

-- When $G$ is compact, the Haar measure is normalized by $m(G)=1$. Then $\widehat{G}$
is discrete, and, if $f\in L^2(G)$, $\chi\in\widehat{G}$, $\widehat{f}(\chi)=\scal{f,\chi}$.
The characters form an orthonormal basis of $L^2(G)$ and $\|\chi\|_\infty=1$.

-- When $G$ is finite, the Haar measure is the counting measure so that $m(G)=|G|$. Then $\widehat{G}\simeq G$
and, if $f\in L^2(G)$, $\chi\in\widehat{G}$, $\widehat{f}(\chi)=|G|^{-1/2}\scal{f,\chi}$.
The characters form an orthogonal basis of $L^2(G)$ and $\|\chi\|_\infty=1$.

\begin{definition}
    Let $\Gamma\subset\widehat{G}$ and $L^2_\Gamma(G)$ be the closure of the span of the functions in $\Gamma$
    in $L^2(G)$-norm. Let $q>2$. The set   $\Gamma$ is called a \emph{$\Lambda_q$-set} if
    $L^q_\Gamma(G)=L^2_\Gamma(G)$, {\it i.e.} there exists a constant $C(q,\Gamma)$ such that, 
    for every $f\in L^2_\Gamma(G)$, $\norm{f}_{L^q(G)}\leq C(q,\Gamma)\norm{f}_{L^2(G)}$.
\end{definition}

As the set of characters $\{\chi\}$ in $\Gamma$ are orthogonal basis for $L^2_\Gamma(G)$, with norm $\norm{\chi}_{L^2(G)}=m(G)^{1/2}$. Therefore,  
every $f\in L^2_\Gamma(G)$ can be written as
$$
f= 
\frac{1}{m(G)^{1/2}}\sum_{\chi\in\Gamma}\widehat{f}(\chi)\chi
$$
with  $\widehat{f}(\chi)= m(G) ^{-1/2}\langle f, \chi\rangle$. Then 
$$
\norm{f}_{L^2(G)}=\left(\sum_{\chi\in\Gamma}|\widehat{f}(\chi)|^2\right)^{\frac{1}{2}}. 
$$

A set $\Gamma\subset\widehat{G}$ is then a $\Lambda_q$ set if, for every $f\in L^2(G)$,
    \begin{equation}
        \label{eq:lambdaqdef}
    \norm{\frac{1}{m(G)^{1/2}}\sum_{\chi\in\Gamma}\widehat{f}(\chi)\chi}_{L^q(G)}\leq C(q,\Gamma)\left(\sum_{\chi\in\Gamma}|\widehat{f}(\chi)|^2\right)^{\frac{1}{2}}.
    \end{equation}  
Of course, when $G$ is finite, every $\Gamma\subset\widehat{G}$ is a $\Lambda_q$-set for any $q>2$.

A celebrated result due to Jean Bourgain \cite{Bourgain89} says the following:
  
\begin{theorem} \label{bourgaintheorem} Let $\Psi=(\psi_1, \dots, \psi_n)$ denote a sequence of $n$ mutually orthogonal functions, with ${\|\psi_i\|}_{L^{\infty}(G)} \leq 1$.   Let $q\geq 2$.
There exists a subset $S$ of $\{1,2, \dots, n\}$, $|S|>n^{\frac{2}{q}}$ such that 
$$ {\left|\left| \sum_{i \in S} a_i \psi_i \right|\right|}_{L^q(G)} \leq \mathcal{B}(q) {\left( \sum_{i \in S} {|a_i|}^2 \right)}^{\frac{1}{2}}.$$ 

The constant $\mathcal{B}(q)$ depends only on $q$ and the estimate above holds for a random set (with respect to the uniform distribution) of size $\lceil n^{\frac{2}{q}} \rceil$, with probability $1-o_N(1)$, where $\lceil x\rceil$ denotes the smallest integer greater than $x$. 
\end{theorem} 

The remarkable feature of this result is that $\mathcal{B}(q)$ depends on $q$ only and not on $n$ nor on $G$.
The evaluation of this constant $\mathcal{B}(q)$ is unfortunately far from obvious from \cite{Bourgain89} and is still
the object of ongoing work.
Specified to characters, this result reads as follows:

\begin{corollary} \label{bourgaindiscretecorollary} 
Let $G$ be a finite group, let 
$q>2$,  and let $q'$ be the dual conjugate, $\dfrac{1}{q'}+\dfrac{1}{q}=1$. Then, with probability $1-o_{|G|}(1)$,
 if $\Sigma$ is a random set of $\widehat{G}$ of size $\left\lceil |G|^{\frac{2}{q}} \right\rceil$,
then
\begin{enumerate}
\renewcommand{\theenumi}{\roman{enumi}}
    \item for every $f\in L^2(G)$ with $\supp \widehat{f}\subset\Sigma$,
$$
\norm{\sum_{\chi\in\Sigma}\widehat{f}(\chi)\chi}_{L^q(G)}
\leq \mathcal{B}_q
\left(\sum_{\chi\in \Sigma}|\widehat{f}(\chi)|^2\right)^{\frac{1}{2}};
$$

    \item for every $f\in L^2(G)$,
    $$
\left(    \sum_{\chi\in\Sigma}|\widehat{f}(\chi)|^2\right)^{1/2}
\leq \mathcal{B}_q|G|^{-1/2}\norm{f}_{L^{q'}(G)}.
$$
\end{enumerate}
\end{corollary}

The first statement states that the synthesis operator 
$T:  \ell^2(\Sigma)\to L^q(G)$
given by 
$T: \displaystyle(a_\chi)_{\chi\in\Sigma}\to \textstyle\sum_{\chi\in\Sigma} a_\chi\chi$
is bounded  with the operator norm $\leq \mathcal{B}_q$. Moreover, the image of this operator is the space of functions whose Fourier support lies in  $\Sigma$. 

The adjoint operator is the analysis map $T^*\,: L^{q'}(G)\to\ell^2(\Sigma)$ given by
$$
T^*f=|G|^{1/2} \bigl(\widehat{f}(\chi)\bigr)_{\chi\in\Sigma}.
$$ 
%
The second statement is therefore just the dual statement of the first one.
It shows that, with high probability, a set of cardinality $\lceil |G|^{\frac{2}{q}} \rceil$
satisfies the $(2,q')$-restriction property with
constant $\rho_{2,q'}(\Sigma)=\mathcal{B}_q|G|$ (independent of $\Sigma$).

We can now reformulate Theorem \ref{th:IMquant}, choosing $\Sigma$ randomly and
taking a set $S$ such that
$$
\mathcal{B}_q|G|^{-\frac{1}{2}}|S|^{\frac{1}{q'}-\frac{1}{2}}=\mathcal{B}_q|G|^{-\frac{1}{2}}|S|^{\frac{1}{2}-\frac{1}{q}}<1.
$$
Writing this more explicitly, we obtain:

\begin{theorem}[Uncertainty Principle in the presence of Randomness]\label{mainbourgain} 
Let $G$ be a finite group, let $q>2$.
Let $\Sigma \subset \widehat{G}$ of size $\left[|G|^{\frac{2}{q}} \right]$, chosen randomly with uniform probability
and $S\subset G$ a set of size
$$
|S|<\mathcal{B}_q^{-\frac{2q}{q-2}}|G|^{\frac{q}{q-2}}.
$$
Then $(S,\Sigma)$ is a strong annihilating pair: for every $f\in L^2(G)$,
$$
\|f\|_{L^2(G)}\leq A_{ann}(S,\Sigma)\bigl(\|f\|_{L^2(G\setminus S)}+\|\widehat{f}\|_{L^2(\widehat{G}\setminus \Sigma)}\bigr).
$$
with
$$
A_{ann}(S,\Sigma)=\frac{1}{1-\mathcal{B}_q|G|^{-\frac{1}{2}}|S|^{\frac{1}{2}-\frac{1}{q}}}
\leq \frac{\mathcal{B}_q+1}{|G|-\mathcal{B}_q^2|S|^{1-\frac{2}{q}}}|G|
$$
and $\mathcal{B}_q$ the constant in Bourgain's theorem.    
\end{theorem}

\begin{remark}
In this theorem, the size condition on the sets $S,\Sigma$ is
$$
|S|^{\frac{q-2}{q}}\,|\Sigma|\lesssim |G|
$$
which is less restrictive than the condition $|S||\Sigma|$ since the $\dfrac{q-2}{q}<1$.
\end{remark}

\begin{remark}
    This result is a quantitative counterpart of previous results by the first and last author of this paper.
    They proved for $G=(\bZ/N\bZ)^d$ that, under the same constraints on the sets $(S,\Sigma)$, they form a weak (thus strong
    annihilating pair). The present proof allows us to compute the annihilation constant.
    This also implies bounds on time-frequency concentration via \eqref{eq:concentration}.
\end{remark}

\subsection{Deterministic examples}

In this section, we will give some explicit examples of our results. We shall need the following notion. 

\begin{definition}
    Let $G$ be a finite abelian group and $E\subset G$. The {\em additive energy} of $E$,
    denoted $\Lambda(E)$ is the number of $(x,y,x',y')\in E^4$ such that $x+y=x'+y'$.
\end{definition}

Recall that $G$ is an additive group while $\widehat{G}$ is multiplicative. So, if $F\subset\widehat{G}$, 
its additive energy is defined by 
$$
\Lambda(F)
=\left|\{\mbox{solutions of }\chi_1\chi_2\overline{\chi_3\chi_4}=1\,:\chi_1,\chi_2,\chi_3,\chi_4\in F\}\right|.
$$

\vskip.125in

The following result is a step in the direction of elucidating the nature of the constant $C_{p,q}$ in the case 
$q=2, p=\dfrac{4}{3}$. To set it up, we need the following result from \cite[Theorem 3.6]{IM24}, which 
 is   adapted here to arbitrary finite groups:

Below, we use the notation $\Bbb Z_N^d = (\Bbb Z/N\Bbb Z)^d$.

\begin{theorem}\label{energyrestrictiontheorem} 
Let $f: {\mathbb Z}_N^d \to {\mathbb C}$ and let $\Sigma$ be a subset of ${\mathbb Z}_N^d$. Then 
\begin{equation} \label{explicitconstantsrestrictionequation} 
\left( \sum_{\chi \in \Sigma} {|\widehat{f}(\chi)|}^2 \right)^{\frac{1}{2}} 
\leq |G|^{-\frac{1}{4}}
\left( \max_{F \subset \Sigma} \frac{\Lambda(F)}{|F|^2} \right)^{\frac{1}{4}} 
\cdot\left( \sum_{x \in {\mathbb Z}_N^d} {|f(x)|}^{\frac{4}{3}} \right)^{\frac{3}{4}}. 
\end{equation} 
In other words,  $\Sigma$ satisfies the $\left(\dfrac{4}{3},2\right)$-restriction estimate with constant
$$
\rho_{4/3,2}(\Sigma)=|G|^{-\frac{1}{4}}
\left( \max_{F \subset \Sigma} \frac{\Lambda(F)}{|F|^2} \right)^{\frac{1}{4}} 
$$
\end{theorem} 

We give the proof for the sake of completeness. 

\begin{proof}
    We write
    $$
\sum_{\chi\in\Sigma}|\widehat{f}(\chi)|^2
=\sum_{\chi\in\widehat{G}}|\widehat{f}(\chi)|^2\mathbf{1}_\Sigma(\chi)
=\sum_{\chi\in\widehat{G}}\widehat{f}(\chi)\mathbf{1}_\Sigma(\chi)g(\chi)
    $$
    with $g(\chi)=\overline{\widehat{f}(\chi)}\mathbf{1}_\Sigma(\chi)$.
By definition
\begin{eqnarray*}
\sum_{\chi\in\widehat{G}}|g(\chi)|^2
   &=&\sum_{\chi\in\Sigma}|\widehat{f}(\chi)|^2\\
  &=&
|G|^{-1/2}\sum_{\chi\in\widehat{G}}\sum_{m\in G}f(m)\overline{\chi(m)}\mathbf{1}_\Sigma(\chi)g(\chi)\\
  &=&\sum_{m\in G}f(m)\left(|G|^{-1/2}\sum_{\chi\in\widehat{G}}\mathbf{1}_\Sigma(\chi)g(\chi)\overline{\chi(m)}\right)\\
  &=&\sum_{m\in G}f(m)\cF_{\widehat{G}}[\mathbf{1}_\Sigma g](m).
\end{eqnarray*}
We then apply H\"older's inequality to obtain
$$
\sum_{\chi\in\Sigma}|\widehat{f}(\chi)|^2
\leq\left(\sum_{m\in G}|f(m)|^{\frac{4}{3}}\right)^{\frac{3}{4}}
\left(\sum_{m\in G}|\cF_{\widehat{G}}[\mathbf{1}_\Sigma g](m)|^4\right)^{\frac{1}{4}}.
$$
It remains to estimate the second factor:
\begin{eqnarray*}
\sum_{m\in G}|\cF_{\widehat{G}}[\mathbf{1}_\Sigma g](m)|^4
&=&\frac{1}{|G|^2}\sum_{m\in G}\sum_{\chi_1,\chi_2,\chi_3,\chi_4\in\Sigma}
g(\chi_1)g(\chi_2)\overline{g(\chi_3)g(\chi_4)} \chi_1(m)\chi_2(m)\overline{\chi_3(m)\chi_4(m)}\\
&=&\sum_{\chi_1,\chi_2,\chi_3,\chi_4\in\Sigma}g(\chi_1)g(\chi_2)\overline{g(\chi_3)g(\chi_4)}
\frac{1}{|G|^2}\sum_{m\in G}\chi_1(m)\chi_2(m)\overline{\chi_3(m)\chi_4(m)}\\
&=&\frac{1}{|G|}\sum_{\chi_1\chi_2=\chi_3\chi_4}g(\chi_1)g(\chi_2)\overline{g(\chi_3)g(\chi_4)}.
\end{eqnarray*}
The modulus of this expression is bounded by
$$
\frac{1}{|G|}\left(\max_{F \subset \Sigma} \frac{\Lambda(F)}{|F|^2}\right)\left(\sum_{\chi\in\widehat{G}}|g(\chi)|^2\right)^2.
$$
In summary, we have shown that
$$
\sum_{\chi\in\widehat{G}}|g(\chi)|^2\leq 
\frac{1}{|G|^{\frac{1}{4}}}\left(\max_{F \subset \Sigma} \frac{\Lambda(F)}{|F|^2}\right)^{\frac{1}{4}}
\left(\sum_{m\in G}|f(m)|^{\frac{4}{3}}\right)^{\frac{3}{4}}
\left(\sum_{\chi\in\widehat{G}}|g(\chi)|^2\right)^{\frac{1}{2}}.
$$
Factoring out the common term on both the left-hand and right-hand sides, we obtain the claimed result
\end{proof}

\vskip.125in 

Using Theorem \ref{energyrestrictiontheorem} and Theorem \ref{th:IMquant}, we obtain the following result. 

\begin{theorem} \label{explicitconstantsannihilatingpairstheorem}   Let $E,S \subset {\mathbb Z}_N^d$ such that 
\begin{equation} \label{universalcondition} 
\max_{U \subset S} \frac{\Lambda(U)}{{|U|}^2} \cdot |E| < N^d.
\end{equation}

Then  for any  $f: {\mathbb Z}_N^d \to {\mathbb C}$ 
\begin{equation} \label{annihilationwithgoodconstants} {||f||}_{L^2({\mathbb Z}_N^d)} \leq C_{ann} \left({||f||}_{L^2(E^c)}{+||\widehat{f}||}_{L^2(S^c)} \right), \end{equation} where $C_{ann}$ may be taken to be 
\begin{equation} \label{goodconstant}  1+\frac{1}{1-\sqrt{\frac{{\left( \max_{U \subset S} \frac{\Lambda(U)}{{|U|}^2} \right)}^{\frac{1}{2}} {|E|}^{\frac{1}{2}}}{N^{\frac{d}{2}}}}}. \end{equation} 

\vskip.125in 

In particular, $(E,S)$ is a strong $L^2$-annihilating pair under the assumptions given above.
\end{theorem} 

Another result that follows readily from our methods is the following $\ell^p({\mathbb Z}_N^d)$ version of Theorem \ref{th:IMquant}. 

\begin{theorem} \label{banachmain} Let $f: {\mathbb Z}_N^d \to {\mathbb C}$. Let $E,S \subset {\mathbb Z}_N^d$ such that $S$ satisfies the $(p,q)$ restriction estimate with norm $C_{p,q}$, for some $1 \leq p \leq 2 \leq q$, and 
\begin{equation} \label{uglycondition} 
{|E|}^{2-p} \cdot |S|<\frac{N^d}{C_{p,q}^p}.\end{equation} 

Then for $1 \leq p \leq 2$, 
\begin{equation} \label{uglymotherfucker} {||f||}_{L^{p'}({\mathbb Z}_N^d)} \leq \frac{N^{-d \left(\frac{1}{2}-\frac{1}{p'} \right)}}{1-{\left( \frac{{|E|}^{2-p}|S|C_{p,q}^p}{N^d}\right)}^{\frac{1}{p}} }{||\widehat{f}||}_{L^p(S^c)}+\left(1+\frac{1}{1-{\left( \frac{{|E|}^{2-p}|S|C_{p,q}^p}{N^d} \right)}^{\frac{1}{p}}}\right){||f||}_{L^{p'}(E^c)}. 
\end{equation}

Since $(1,q)$ restriction estimate always holds with $C_{1,q}=1$, then for any sets $E,S \subset {\mathbb Z}_N^d$ such that $|E||S|<N^d$, 
\begin{equation} \label{dailyuglymotherfucker} {||f||}_{L^{\infty}({\mathbb Z}_N^d)} \leq \frac{N^{-\frac{d}{2}}}{1-\frac{|E||S|}{N^d}} {||\widehat{f}||}_{L^1(S^c)}+\left(1+\frac{1}{1-\frac{|E||S|}{N^d}} \right) {||f||}_{L^{\infty}(E^c)}. \end{equation} 
\end{theorem} 

\vskip.125in  




\section{Euclidean examples: $G=\R^d$.}\label{Sec:excont}


Suppose $\bS$ is a smooth hypersurface of $\R^n$, 
$\mbox{d}\sigma_{\bS}$ its surface measure and that the $p,q$-restriction estimate holds on $\bS$ in the sense that there exists a constant $C_{p,q}(\bS)$ such that, for every $f\in\mathcal{S}(\R^d)$
\begin{equation}
    \label{eq:resconj}
\left(\int_{\bS}|\widehat{f}(\xi)|^q\,\mbox{d}\sigma_{\bS}(\xi)\right)^{\frac{1}{q}}
\leq C_{p,q}(\bS)\left(\int_{\bR^d}|f(x)|^p\,\mbox{d}x\right)^{\frac{1}{p}}.
\end{equation}

\vskip.125in 

Let $\bS_\delta=\{\xi\in\R^d\,:\ \exists\eta\in\bS,\ |\eta-\xi|<\delta/2\}$, the $\delta$-neighborhood of $\bS$. Then, for $\delta$ small enough, $\bS$ satisfies the $p,q$-restriction with constant $\rho_{p,q}(\bS)=C_{p,q}(\bS)\delta^{1/q}$, in the sense that for every $f\in\mathcal{S}(\R^d)$,
$$
\left(\int_{\bS_\delta}|\widehat{f}(\xi)|^q\,\mbox{d}\xi\right)^{\frac{1}{q}}
\leq C_{p,q}(\bS)\delta^{\frac{1}{q}}\left(\int_{\bR^d}|f(x)|^p\,\mbox{d}x\right)^{\frac{1}{p}}.
$$

We can scale the problem further. If $R>0$ then 
$$(R\bS)_\delta= R(\bS_{\delta/R}) \ \text{and} \ \widehat{f}(R\xi)=\widehat{f_R}(\xi),$$ with $f_R(x)=R^{-d}f(x/R).$ It follows that 
\begin{eqnarray*}
\left(\int_{(R\bS)_{\delta/R}}|\widehat{f}(\xi)|^q\,\mbox{d}\xi\right)^{\frac{1}{q}}
&=&R^{\frac{d}{q}}\left(\int_{\bS_{\delta/R}}|\widehat{f}(R\eta)|^q\,\mbox{d}\eta\right)^{\frac{1}{q}}\\
&\leq&R^{\frac{d}{q}}C_{p,q}(\bS)\left(\frac{\delta}{R}\right)^{\frac{1}{q}}
\left(\int_{\R^d}|R^{-d}f(x/R)|^p\,\mbox{d}x\right)^{\frac{1}{p}}\\
&=&R^{d\left(\frac{1}{p}+\frac{1}{q}-1\right)}\left(\frac{\delta}{R}\right)^{\frac{1}{q}}C_{p,q}(\bS)\left(\int_{\R^d}|f(y)|^p\,\mbox{d}x\right)^{\frac{1}{p}}.
\end{eqnarray*}
It follows that $(R\bS)_\delta$ satisfies the $(p,q)$-restriction estimate
with 
$$
\rho_{p,q}\bigl((R\bS)_\delta\bigr)=C_{p,q}(\bS)
R^{d\left(\frac{1}{p}+\frac{1}{q}-1\right)}\left(\frac{\delta}{R}\right)^{\frac{1}{q}}.
$$
Notice that, if $p<q$, then $\rho_{p,q}\bigl((R\bS)_\delta\bigr)\to 0$ as $\delta\to 0$. 

\vskip.125in 

Recall that the restriction conjecture for the sphere states the following:
\begin{conjecture}
When $\bS=\bS^{d-1}$ is the sphere on $\R^d$, then \eqref{eq:resconj} is conjectured to hold
when $p,q$ satisfy
$$
p<\dfrac{2d}{d+1}\quad\mbox{and}\quad q\leq \dfrac{d-1}{d+1}p'
$$
where as usual $\dfrac{1}{p}+\dfrac{1}{p'}=1$.
\end{conjecture}

Despite its resolution in dimension $d=2$ in work by Fefferman \cite{F69} and Zygmund \cite{Z74}, the conjecture remains open in dimension $d\geq 3$. We shall need the following result established by Tomas, with the endpoint obtained by Stein. 

\begin{theorem}[Stein-Tomas] When $\bS=\bS^{d-1}$ is the sphere on $\R^d$, then \eqref{eq:resconj} holds for 
$$
1\leq p\leq p_{TS}:=\dfrac{2(d+1)}{d+3} \ \text{and} \ q=2.
$$
\end{theorem}

Now, let 
$$
\bA(R,\delta)=(R\bS^{d-1})_\delta=\left\{x\in\R^d\,: R-\dfrac{\delta}{2}<|x|<R+\dfrac{\delta}{2}\right\}
$$
an annulus and notice that, for $R$ large enough and $\delta$ small enough,
$$
|\bA(R,\delta)|\approx R^{d-1}\delta
$$
and
$$
\rho_{p_{TS},2}\bigl(\bA(R,\delta)\bigr)
=C_{p_{TS},2}(\bS)R^{\frac{1}{2}\frac{d-1}{d+1}}\delta^{\frac{1}{2}}
=\kappa_d|\bA(R,\delta)|^{\frac{1}{2(d+1)}}\delta^{\frac{d}{2(d+1)}}
$$
Here $C_{p_{TS},2}(\bS)$ is the restriction 
constant given by the Stein-Tomas Theorem and thus $\kappa_d$
is a constant that depends on the dimension $d$ only.

\vskip.125in 

If $S$ is a set of finite measure, then 
according to Nazarov's 
uncertainty principle in higher dimension due to the second author, $S,\bA(R,\delta)$ form a strong annihilating pair with annihilating constant
$$
C\bigl(S,\bA(R,\delta)\bigr)=ce^{c|S||\bA(R,\delta)|}
\asymp
ce^{c'|S|R^{d-1}\delta}.
$$
Note that it is conjectured that this can be improved to $ce^{c(|S||\bA(R,\delta)|)^{1/d}}$.

Let us show that if $\delta$ is sufficiently small compared to $R$, one can improve this with
Theorem \ref{th:IMquant}. To apply this theorem, we need
$$
\rho_{p_{TS},2}\bigl(\bA(R,\delta)\bigr)|S|^{\frac{1}{p_{TS}}-\frac{1}{2}}<1
$$
that is
$$
|\bA(R,\delta)|^{1/2}|S|\delta^{\frac{d}{2}}\leq
\kappa_d R^{\frac{d-1}{2}}\delta^{\frac{d+1}{2}}|S|<1
$$
where $\kappa_d$ is a constant that depends on the dimension $d$ only. Theorem \ref{th:IMquant}
then shows that
$$
C\bigl(S,\bA(R,\delta)\bigr)=1+\frac{1}{1-\kappa_d R^{\frac{d-1}{2}}\delta^{\frac{d+1}{2}}|S|}.
$$

Now fix $S$ with $|S|>0$, $R>0$ and $\delta=(2\kappa_d|S|)^{-\frac{2}{d+1}}R^{-\frac{d-1}{d+1}}$ then 
$$
|\bA(R,\delta)|\approx \left(\frac{R^{d(d-1)}}{|S|^2}\right)^{\frac{1}{d+1}},\ 
|S||\bA(R,\delta)|\approx R^{(d-1)\frac{d}{d+1}}|S|^{\frac{d-1}{d+1}}
\ \mbox{while}\ 
C\bigl(S,\bA(R,\delta)\bigr)=3.
$$
So the measure of $S$ and of $\bA(R,\delta)$ can be arbitrarily large but
the annihilation constant stays constant, in strong contrast with the constant obtained
in Nazarov's theorem. In summary, we obtain the following result: 

\begin{proposition}
Let $d\geq 2$, let $S\subset\R^d$ be a set of positive measure, $R>0$ and $\delta>0$.
Then, if $R\gtrsim 1$ and $d\lesssim |S|^{-\frac{2}{d+1}}R^{-\frac{d-1}{d+1}}$,
$\bigl(S,\bA(R,\delta)\bigr)$ is a strong annihilating pair with annihilation constant $3$. In particular, for every $f\in L^2(\R^d)$,
$$
\norm{f}_{L^2(\R^d)}\leq 3\bigl(\norm{f}_{L^2(\R^d\setminus S)}+
\norm{f}_{L^2(\R^d\setminus \bA(R,\delta))}\bigr).
$$
\end{proposition}


\section*{Acknowledgements}

The first listed author was supported in part by the National Science Foundation DMS - 2154232. The third listed author was supported in part by AMS-Simons Research Enhancement Grant, Simons Foundation Fellowship, and the PSC-CUNY research grants. The second listed author was supported in part by ANR24-CE40-5470 grant.

The authors would like to thank the Isaac Newton Institute for Mathematical Sciences, Cambridge, for support and hospitality during the program entitled ``Multivariate approximation, discretization, and sampling recovery," where work on this paper was undertaken. This work was supported by EPSRC grant EP/R014604/1.

 \end{document}